\documentclass[12pt,a4paper]{amsart}
\usepackage[utf8]{inputenc}
\usepackage{amsmath}
\usepackage{amsfonts}
\usepackage{amssymb}
\usepackage{color}
\usepackage{cite}
\usepackage[english]{babel}

\newtheorem{theorem}{Theorem}
\newtheorem{lemma}{Lemma}

\newcommand{\Rd}{{\mathbb R}^d}
\newcommand{\RR}{\mathbb R}
\newcommand{\NN}{\mathbb N}

\DeclareMathOperator{\mes}{meas}
\DeclareMathOperator{\diam}{diam}

\allowdisplaybreaks
\begin{document}

\title[On optimization of cubature formulae]{On optimization of cubature formulae for Sobolev classes of functions defined on star domains}
\author[O.~V.~Kovalenko]{Oleg Kovalenko}
\address{Department of Mathematical Analysis and Theory of Functions, Oles Honchar Dnipro National University, Dnipro, Ukraine}
\email{olegkovalenko90@gmail.com}

\begin{abstract}
We find asymptotically optimal methods of recovery of the integration operator given values of the function at a finite number of points for a class of multivariate functions defined on a bounded star domain that have bounded in $L_p$ norm of their distributional gradient; thus we generalize the known solution of this optimization problem in the case, when  the domain of definition of the functions is convex.  
\end{abstract}
\maketitle

\section{Introduction}
Let a bounded measurable set $Q\subset\mathbb{R}^d$, a class $X$ of continuous on $Q$ functions, and $n\in\mathbb{N}$ be given. Arbitrary function $\Phi\colon \mathbb{R}^n\to \mathbb{R}$ is called a method of recovery. For given points $x_1,\dots,x_n\in Q$ the error of recovery of the integral by the method $\Phi$ is
$$e(X, \Phi, x_1,\dots,x_n):= \sup_{f\in X}\left|\,\int_{Q}f(x)dx - \Phi(f(x_1),\dots,f(x_n))\right|.$$

The problem of the optimal recovery of the integral is to find the best error of recovery 
\begin{equation}\label{optimalError}
E_n(X):=\inf_{x_1,\dots,x_n\in Q}\inf_{\Phi\colon\mathbb{R}^n\to\mathbb{R}} e(X,\Phi,x_1,\dots,x_n),
\end{equation}
the best method of recovery, and the best position of the informational set $x_1,\dots,x_n$ i.e., such method $\tilde{\Phi}\colon \mathbb{R}^n\to \mathbb{R}$ and points $\tilde{x}_1,\dots,\tilde{x}_n\in Q$, for which the infima in~\eqref{optimalError} are attained (if such a method and points exist). 

In many cases it is hard to find the best error of recovery and an optimal recovery method; in such situations it is interesting to find an asymptotically optimal method of recovery i.e., such sequence of methods $\Phi_n\colon\RR^n\to\RR$ and informational sets $\{x_1^n,\dots, x_n^n\}$,  $n\in\NN$, that 
$$\lim_{n\to\infty}\frac{E_n(X)}{e(X, \Phi_n, x_1^n,\dots,x_n^n)} = 1.$$

The problem of optimal recovery and, in particular, the problem of optimization of cubature formulae has a rich history, see e.g. monographs~\cite{Traub80, Zhensykbaev,Plaskota,Osipenko2000}.

Let $Q\subset \Rd$, $d\in\NN$, be a nonempty bounded open set. By  $W^{1,p}(Q)$, $p\in [1,\infty]$, we denote the Sobolev space of functions $f\colon Q\to \RR$ such that $f$ and all their (distributional) partial derivatives of the first order belong to $L_p(Q)$. 

As usually, for $x=(x^1,\dots, x^d)\in \Rd$ and $q\in [1,\infty)$ set $|x|_q:= \left(\sum_{k=1}^d|x^k|^q\right)^\frac {1}{q}$, $|x|_\infty:= \max_{k=1,\dots, d}|x^k|$. It is clear that for all $f\in W^{1,p}(Q)$ we have $\|\,|\nabla f|_1\,\|_{L_p(Q)}<\infty$. For $p\in[1,\infty]$ set $$W^{\infty}_{p}(Q):=\{f\in W^{1,p}(Q)\colon \|\,|\nabla f|_1\,\|_{L_p(Q)}\leq 1\}.$$

Below we consider only the case $p>d\geq 2$. For $p>d$ imbedding of the class $W^{1,p}(Q)$ into the space of bounded continuous on $Q$ functions holds, provided by some restrictions on the geometry of $Q$ (sets $Q$ for which the imbedding holds will be called admissible). For example, it is sufficient to require that $Q$ satisfies the cone condition, see Chapter 4 and Theorem~4.12 in~\cite{Adams}.

A bounded open set $Q\subset\Rd$ is called a star domain with respect to a ball $B$ (or simply a star domain for brevity), if for all $x\in Q$ and $y\in B$ the segment $xy$ belongs to $Q$. It is not hard to verify that the interior of the closure of a star domain $Q$ coincides with $Q$. This implies that the closure of $Q$ is an asymmetric star body according to~\cite[Definition~I.2.2]{geometryOfNumbers}, hence its distance function is continuous (see~\cite[Theorem~I.2.2]{geometryOfNumbers}), and thus $Q$ is Jordan measurable (see the proof of~\cite[Theorem~I.1.5]{geometryOfNumbers}). 
Moreover, it is easy to see that a star domain satisfies the cone condition, and hence each function $f$ from $W_p^\infty(Q)$ has a continuous representation.

Everywhere below, for a finite set $A$, by $|A|$ we denote the number of its elements. For $x,y\in \Rd$ by $(x,y)$ we denote the dot product of $x$ and $y$. For a bounded set $Q\subset \Rd$ we set $\diam Q:= \sup_{a,b\in Q}|a - b|_2$. We write $\mes Q$ to denote the Lebesgue measure of a measurable set $Q$.

The goal of the article is to find asymptotically optimal cubature formulae on the class $W^\infty_p(Q)$, where $Q$ is a bounded star domain. In the case, when $Q$ is a convex bounded domain, such problem was solved in~\cite{Kovalenko20b}. Solutions to some extremal problems for similar classes of functions can be found in~\cite{BP2012-nabla,Kovalenko21b,Kovalenko22a}.

\section{The main result}
\subsection{Asymptotically optimal informational sets and methods}
The construction of asymptotically optimal informational sets and recovery methods in the case when $Q$ is a star domain is the same as in the case of a convex set $Q$. We adduce it below together with some properties that we will need, see~\cite{Babenko76, Kovalenko20b,Kovalenko21b}.

For a given $h>0$ consider the lattice $\Lambda$ in $\Rd$ generated by the vectors $(2h,0,0,\dots,0)$, $(0,2h,0,0,\dots,0),\dots,$ $(0,\dots,0,2h)\in\Rd$. Denote by $P_k(h)$, $k\in\mathbb{N}$, the cubes into which the lattice $\Lambda$ divides $\Rd$; their volumes are equal to $(2h)^d$. Denote by $A(h)$ the set of all cubes $P_k(h)$ that are contained in $Q$, let $a(h)$ be the set of the centers of the cubes from $A(h)$. Denote by $B(h)$ the set of all cubes $P_k(h)$ that have non-empty set of common with $Q$ internal points.
For each $n\in\mathbb{N}$ set 
\begin{equation}\label{h_n}
h_n:=\frac 12 \left(\frac{\mes  Q}{n}\right)^{\frac 1d}.
\end{equation}

For each cube $P$ from the set $B(3h_n)$ choose a point by the following rule: the center of the cube $P$, if it belongs to $a(h_n)$; else, arbitrary point from $P\cap a(h_n)$, if the intersection is not empty; else, arbitrary internal point of $Q\cap P$.
  
  Denote by $S_1(n)$ the set of chosen points;
by $S_2(n)$ arbitrary subset of the set $a(h_n)\setminus S_1(n)$ that contains $n-|S_1(n)|$ points (for large enough $n$ this number is positive; if the set $a(h_n)\setminus S_1(n)$ contains less than  $n-|S_1(n)|$ points, then we take all points of $a(h_n)\setminus S_1(n)$). Set 
\begin{equation}\label{S_n}
S(n):=S_1(n)\cup S_2(n).
\end{equation}

Let $S(n)=\{x_1^*,\dots, x_{|S(n)|}^*\}.$ For each $k=1,\dots, |S(n)|$ define the set
\begin{equation}\label{Qpartition}
V_k:=\{x\in Q\cap P(3h_n;x_k^*)\colon |x-x_k^*|_\infty < |x-x_s^*|_\infty, s\neq k\},
\end{equation}
where $P(3h_n;x_k^*)$ is the cube from $B(3h_n)$ that contains $x_k^*$.
 Then the sets $V_k$ are pairwise disjoint, $\bigcup_{k=1}^{|S(n)|}V_k\subset Q$ and $\mes \left(Q\setminus\bigcup_{k=1}^{|S(n)|}V_k\right) = 0$. 
 Set 
\begin{equation}\label{optimalWeights}
c_k^*:=\mes  V_k, k=1,\dots, |S(n)|.
\end{equation}

The following lemma states some of the properties of the sets $S(n)$ and $V_k$, defined above, see~\cite[Lemma~4]{Kovalenko20b}.

\begin{lemma}\label{l::InformationalSetProps}
	Let $Q\subset \Rd$ be a Jordan measurable set, $n\in\mathbb{N}$ be enough big and $h_n$ be defined by~\eqref{h_n}. Then the following properties hold:
	\begin{enumerate}
		\item $S(n)\subset Q$ and $|S(n)|\leq n$.
		\item\label{VkDiam} If $x\in V_k$ then $|x-x_k^*|_\infty \leq 6h_n$.
		\item\label{OnePointInCube} For each cube $P\in B(h_n)$, $|P\cap S(n)|\leq 1$.
		\item\label{U_k} Let $R_n$ be the union of cubes $P\in A(h_n)$ with centers that belong to $S(n)$. Denote by $U_k:=V_k\setminus R_n$. Then $\mes \bigcup_{k=1}^{|S(n)|}U_k = o(1)$, $n\to\infty$.
	\end{enumerate}
\end{lemma}

 \subsection{Asymptotically optimal cubature formulae}
The following theorem is the main result of the article. 
\begin{theorem}\label{th::ORunitWeight}
	Let $d\geq 2$, $p\in(d,\infty]$ and a bounded star domain $Q$ be given. Then 
	 $$E_n\left(W_{p}^{\infty}(Q)\right)=c(d,p)\left(\frac {\mes Q}{2^d}\right)^{\frac 1 d +\frac 1 {p'}}\cdot \frac{1+o(1)} {n^{\frac 1 d}},\, n\to\infty,$$
where $c(d,p):=\frac 1 d \left\|\frac {1} {|\cdot|^{d-1}_\infty} - |\cdot|_\infty\right\|_{L_{p'}(\{x\in\Rd\colon |x|_\infty\leq 1\})}$. 

The asymptotically optimal informational set is $S(n)$ defined by~\eqref{S_n}, and the optimal recovery method is 
\begin{equation*}
	 \tilde{\Phi}_n(f(x_1),\dots, f(x_{|S(n)|}))=\sum_{k=1}^{|S(n)|} c_k^*f(x_k),
\end{equation*}
  where the weights $c_k^*$ are defined by~\eqref{optimalWeights}.
\end{theorem}

It was proved, see~\cite[Lemma~5]{Kovalenko20b} that for arbitrary admissible $Q$
$$E_n\left(W_{p}^{\infty}(Q)\right)
\geq 
c(d,p)\left(\frac {\mes Q}{2^d}\right)^{\frac 1 d +\frac 1 {p'}}\cdot \frac{1+o(1)} {n^{\frac 1 d}},\, n\to\infty$$
and, see~\cite[Lemma~6]{Kovalenko20b}, as $n\to\infty$,
$$
E_n\left(W_{p}^{\infty}(Q)\right)
\leq 
\sup_{f\in W_{p}^{\infty}(Q)}\left|\, \sum_{k=1}^{|S(n)|}\int_{U_k} [f(x)-f(x_k^*)]dx \right| 
+ 
c(d,p)\left(\frac {\mes Q}{2^d}\right)^{\frac 1 d +\frac 1 {p'}}\cdot \frac{1+o(1)} {n^{\frac 1 d}},
$$
where the sets $U_k$ are defined in Property~\ref{U_k} of Lemma~\ref{l::InformationalSetProps}.
Thus in order to prove the theorem it is sufficient to prove the following lemma.
\begin{lemma}\label{l::th5RemainderEstimate}
Let $d\in\mathbb{N}$, $p\in(d,\infty]$ and a bounded star domains $Q$ be given. Then
\begin{equation}\label{th6.remainderEstimate}
\sup_{f\in W_{p}^{\infty}(Q)}\sum_{k=1}^{|S(n)|}\int_{U_k} \left|f(x)-f(x_k^*)\right|dx = o\left(n^{-\frac 1 d}\right),\;n\to\infty.
\end{equation}
\end{lemma} 

A crucial tool in the proof of the main results is the following result, see~\cite[Chapter 6.9]{Lieb}.
\begin{lemma}\label{l::LagrangeTh}
Suppose $p>d$ and $Q\subset \Rd$ is admissible. Let $f\in W^{1,p}(Q)$ and $x,y\in Q$ be such that the whole segment with ends in the points $x$ and $y$ belongs to $Q$. Then
$$f(y)-f(x)=\int_0^1\left(y-x,\nabla f[(1-t)x+ty]\right)dt.$$
\end{lemma}
Observe that since the sets $U_k$, $k=1,\ldots, |S(n)|$ are generally speaking not convex, we can't directly apply Lemma~\ref{l::LagrangeTh} to the difference $f(x) - f(x_k^*)$ under the integral in~\eqref{th6.remainderEstimate}.
We define functions $p_k\colon U_k\to Q$, $k=1,\dots, |S(n)|$, such that whole segments $xp_k(x)$ and $p_k(x)x_k^*$ belong to $Q$ and they are ''not much longer'' than the segment $xx_k^*$. Once this is done we can write the inequality 
$$|f(x)-f(x_k^*)|
\leq 
|f(x)-f(p_k(x))|
+
|f(p_k(x))-f(x_k^*)|,
$$
and apply Lemma~\ref{l::LagrangeTh} to switch from the values of the function $f$ to the values of its gradient. This will allow to obtain an estimate from above for the quantity~\eqref{th6.remainderEstimate} in terms of $\|\,|\nabla f|_1\,\|_{L_p(Q)}$ and the total measure of the sets $U_k$, which in turn will imply Lemma~\ref{l::th5RemainderEstimate}.
\section{Proof of the main result}
\subsection{Auxiliary results}

We need the following lemmas. 
\begin{lemma}\label{l::integralsSum}
Let $T_k\subset U_k$ and functions $\phi_k\colon T_k\to Q$  be such that $\phi_k(T_k)$ is measurable, $k=1,\dots, |S(n)|$. Assume that there exists a number $c>0$ such that $|\phi_k(x) - x_k^*|_\infty\leq c|x - x_k^*|_\infty$ for all $x\in T_k$, $k=1,\dots |S(n)|$. Then there exists a number $C>0$ that does not depend on $n$ and such that for all integrable on $Q$ functions $g$
$$\sum_{k=1}^{|S(n)|}\int_{\phi_k(T_k)}|g(x)|dx \leq C\int_Q|g(x)|dx.$$
\end{lemma}

\begin{proof}
If the numbers $1\leq k_1<k_2\leq |S(n)|$ and the points $x\in T_{k_1}$ and $y\in T_{k_2}$ are such that $z=\phi_{k_1}(x) = \phi_{k_2}(y)$, then, using Property~\ref{VkDiam} of Lemma~\ref{l::InformationalSetProps}, we obtain
\begin{multline}\label{centersDist}
|x_{k_1}^*-x_{k_2}^*|_\infty\leq |x_{k_1}^*-z|_\infty + |z-x_{k_2}^*|_\infty = |x_{k_1}^*-\phi_{k_1}(x)|_\infty 
 + |\phi_{k_2}(y)-x_{k_2}^*|_\infty 
\\ \leq 
c|x_{k_1}^*-x|_\infty + c|y- x_{k_2}^*|_\infty \leq 12 h_n\cdot c.
\end{multline}

It is now sufficient to prove that there exists a number $N\in\mathbb{N}$ that does not depend on $n$, and a partition of the set $\{T_1,\dots, T_{|S(n)|}\}$ into groups $\{T_{k_1^i},\dots, T_{k_{m_i}^i}\}$, $i=1,\dots, N$, such that for all $i=1,\dots, N$ and different $j,s\in\{k_1^i,\dots, k_{m_i}^i\}$, $|x_j^*-x_s^*|_\infty > 12h_n\cdot c$. Really, if such partition is done, then, due to~\eqref{centersDist}, the sets  $\{\phi_{k_1^i}(T_{k_1^i}),\dots, \phi_{k_{m_i}^i}(T_{k_{m_i}^i})\}$ are pairwise disjoint for each $i = 1\dots, N$. Hence we obtain
\begin{multline*}
\sum_{k=1}^{|S(n)|}\int_{\phi_k(T_k)}|g(x)|dx 
= \sum_{i=1}^{N}\sum_{s=1}^{m_i}\int_{\phi_{k_s^i}(T_{k_s^i})}|g(x)|dx
\\ = 
\sum_{i=1}^{N}\int_{\bigcup_{s=1}^{m_i}\phi_{k_s^i}(T_{k_s^i})}|g(x)|dx
\leq \sum_{i=1}^{N}\int_{Q}|g(x)|dx = N\int_{Q}|g(x)|dx. 
\end{multline*}

To do such partition we associate an index $I_P\in\mathbb{Z}^d$ with each of the cubes $P\in B(h_n)$. The index $I_P$ is equal to the coordinates of the ''left bottom point'' of $P$ (i.e., the point of the cube $P$ with minimal coordinates) in the basis of the lattice. Let $M$ be an integer bigger than $6c + 1$. We divide all cubes from $B(h_n)$ into $N= M^d$ groups in such a way that two cubes $P_1,P_2\in B(h_n)$ belong to the same group if and only if all the coordinates of $I_{P_1}-I_{P_2}$ are divisible by $M$.

Let two different cubes $P_1,P_2$ belong to one group and $x\in P_1$, $y\in P_2$. Then due to the definition of the group $|x-y|_\infty\geq (M-1)\cdot 2h_n>12ch_n$. 

Now we construct a partition of the sets $T_k$, $k=1,\dots, |S(n)|$. We put two sets $T_k$ and $T_j$ into one group if and only if $x_k^*$ and $x_j^*$ belong to cubes from the same group. Such partition is a desired one, since Property~\ref{OnePointInCube} of Lemma~\ref{l::InformationalSetProps} holds. The lemma is proved.\end{proof}
\begin{lemma}\label{l::det}
	Let two vectors $u=(u^1,\dots, u^d)$, $v=(v^1,\dots, v^d)$ and numbers $\alpha, \beta\in\mathbb{R}$ be given. Then 
	$$\det\left(\alpha\cdot I + \beta\cdot \left\|u^iv^j\right\|_{i,j=1}^d\right) = \alpha^{d-1}(\alpha + \beta (u,v)),$$
	where I denotes the identity matrix.
\end{lemma}

\begin{proof} We prove the lemma by induction on $d$. If $d=2$, then 
$$\begin{vmatrix}
  \alpha + \beta u^1v^1 & \beta u^1v^2 \\
  \beta u^2v^1  & \alpha + \beta u^2v^2
\end{vmatrix}
=
\alpha^2 + \alpha\beta u^1v^1+\alpha\beta u^2v^2=\alpha(\alpha + \beta (u,v)).$$

Let the lemma be true for some $d-1\in \mathbb{N}$. Then
\begin{gather*}
\det\left(\alpha\cdot I + \beta\cdot \left\|u^iv^j\right\|_{i,j=1}^d\right)
\\=
\begin{vmatrix}
  \alpha  & 0 & \cdots  & 0 \\
  \beta u^2v^1  & \alpha + \beta u^2v^2 & \cdots & \ \beta u^2v^d \\
  \vdots & \vdots & \ddots & \vdots \\
  \beta u^dv^1 & \beta u^dv^{2} & \cdots &   \alpha + \beta u^dv^d
\end{vmatrix}
+
\beta u^1
\begin{vmatrix}
  v^1 & v^2 & \cdots  & v^d \\
  \beta u^2v^1  & \alpha + \beta u^2v^2 & \cdots  & \beta u^2v^d \\
  \vdots & \vdots & \ddots &\vdots  \\
  \beta u^dv^1 & \beta u^dv^{2} & \cdots  &  \alpha + \beta u^dv^d
\end{vmatrix}
\\=
\alpha\cdot \alpha^{d-2}\left(\alpha + \beta \sum_{k=2}^{d}u^kv^k\right) + \beta u^1
\begin{vmatrix}
  v^1 & v^2 & \cdots & v^d \\
  0  & \alpha  & \cdots & 0 \\
  \vdots & \vdots & \ddots & \vdots \\
  0 & 0 & \cdots &  \alpha 
\end{vmatrix}
\\=
\alpha^{d-1}\left(\alpha + \beta \sum_{k=2}^{d}u^kv^k\right) + \beta u^1\alpha^{d-1}v^1
=\alpha^{d-1}(\alpha + \beta (u,v)).
\end{gather*}
The lemma is proved. \end{proof}

\subsection{Auxiliary construction}
Let $Q$ be a star domain with respect to a ball $S_R^d(o)$ with radius $R>0$ and center $o\in Q$. For all  $k=1,\dots, |S(n)|$ and $0<r<R$ define  a function $p_k(\cdot;r)\colon U_k\to\Rd$ by the following equation
\begin{equation}\label{p_kDefinition}
p_k(x;r):=\frac{r\cdot(x+x_k^*)+|x-x_k^*|_2\cdot o}{|x-x_k^*|_2 + 2r}.
\end{equation}

Everywhere below we assume that the distance from $o$ to the boundary $\partial Q$ of the set $Q$ is greater than $R$ (otherwise we can decrease $R$). We also consider so large $n\in\mathbb{N}$ that all the points $x_k^*$ (from $V_k$ with non-empty $U_k$) and sets $U_k$ are outside  of the ball $S_R^d(o)$. As the value of $r$ we will use either $r = R$, or $r = \frac R8$, so that $r$ will be separated from $0$ and independent of $n$.

Equality~\eqref{p_kDefinition} has the following geometrical sense. 

\begin{lemma}\label{pr::geometricSense}
Assume that a point $x\in U_k$ is such that vectors $\overline{ox}$ and $\overline{ox_k^*}$ are not collinear. Consider the $2$-dimensional space $E^2$ generated by these two vectors. Let $S^2_r(o)$ be the boundary circle of  $E^2\cap S_r^d(o)$. Let $o_1o_2$ be the diameter of the circle $S^2_r(o)$ parallel to the segment $xx_k^*$. Then $p_k(x;r)$ is the point where the diagonals of the convex hull of the points $o_1$, $o_2$, $x$ and $x_k^*$ intersect.
\end{lemma}
This lemma will be proved below together with the following lemma.
\begin{lemma}\label{pr::smallDistance}
Segments $xp_k(x)$ and $p_k(x)x_k^*$ are fully contained in $Q$ and there exists a number $C$ that does not depend on $n$ such that 
$$\max \{|x-p_k(x;r)|_\infty, |x_k^*-p_k(x;r)|_\infty\}\leq C {|x-x_k^*|_\infty}.$$
\end{lemma}

\begin{proof}
Let $p$ be the intersection of the diagonals $xo_1$ and $x_k^*o_2$. Then triangles $xpx_k^*$ and $o_1po_2$ are similar. This, in particular, means that the points $o$, $p$ and the middle $m$ of the segment $xx^*_k$ belong to a line and 
$\frac{|m-p|_2}{|o-p|_2} = \frac{|x-x_k^*|_2}{2r}.$
Hence 
$$p=\frac 1 {2r+|x-x_k^*|_2}(|x-x_k^*|_2\cdot o + 2r\cdot m)=\frac{r\cdot(x+x_k^*)+|x-x_k^*|_2\cdot o}{|x-x_k^*|_2 + 2r}.
$$
This means that $p=p_k(x;r)$ and Lemma~\ref{pr::geometricSense} is proved.

From the similarity of triangles $xpx_k^*$ and $o_1po_2$ it also follows that 
$\frac{|x-p|_2}{|o_1-p|_2} =  \frac{|x-x_k^*|_2}{2r}.$
Hence $$|x-p|_\infty\leq |x-p|_2=\frac{|o_1-p|_2\cdot |x-x_k^*|_2}{2r}\leq \frac{\diam Q\cdot \sqrt{d}\cdot |x-x_k^*|_\infty}{2r}.$$
Hence the inequality in the statement of the lemma holds with the constant $C:=\frac {\diam Q \sqrt{d}} {2r}$. The estimate for $|x_k^*-p_k(x,r)|_\infty$ can be obtained using the same arguments. 

Segments $xp_k(x)$ and $p_k(x)x_k^*$ are fully inside $Q$, since the set $Q$ is a star domain and they are subsegments of the segments $xo_1$ and $x_k^*o_2$. Lemma~\ref{pr::smallDistance} is proved.
\end{proof}
\subsection{Some properties of the functions $p_k(\cdot;r)$}
Below we state several properties of the functions  $p_k(\cdot;r)$ defined in the previous paragraph. Everywhere in this paragraph we assume $k\in\{1,\dots, |S(n)|\}$ and $0<r\leq R$ are fixed, and $n$ is big enough, so that the sets $U_k$ are outside the ball $S_R^d(o)$.

\begin{lemma}\label{pr::phiPsiInjectiveness}
For each $t\in [0,1]$ there exists a partition $U_k = \bigcup_{i=1}^4 U_k^i(t)$ into measurable sets $U_k^i(t)$, $i =1,\ldots, 4$, such that the function
\begin{equation}\label{psiDefinition}
    \psi_k(x;r,t):=t\cdot x + (1-t)\cdot p_k(x;r)
\end{equation} is injective on each of the sets $U_k^1(t),\ldots, U_k^4(t)$. For all $t\in(0,1]$ the functions
\begin{equation}\label{phiDefinition}
    \varphi_k(x;r,t):=(1-t)\cdot x_k^* + t\cdot p_k(x;r),
\end{equation}
are injective on each of the sets $U_k^1(0),\ldots, U_k^4(0)$.
\end{lemma}
\begin{proof}
First of all note that for $x,y\in Q$ and $t\in (0,1]$, 
$$
    \varphi_k(x;r,t) = \varphi_k(y;r,t)\iff p_k(x;r) = p_k(y;r)
    \iff 
    \psi_k(x;r,0) = \psi_k(y;r,0).
$$
Thus the  statement about functions~\eqref{phiDefinition} follows from the statement about  functions~\eqref{psiDefinition}.

Let $t\in[0,1]$ and $q\in \psi_k(U_k;r,t)$ be fixed. Next we prove that the preimage $\psi_k^{-1}(q;r,t)$ consists of at most $4$ points.

If the points $q$, $o$ and $x_k^*$ belong to one line (i.e., are linearly dependent), then all the points from the set $\psi_k^{-1}(q;r,t)$ also belong to this line and the proof of the lemma is analogous to the proof below for the case when the points $q$, $o$ and $x_k^*$ are linearly independent.

Assume that the points $q$, $o$ and $x_k^*$ are linearly independent. Then they generate a 2-dimensional space. From the geometrical sense of $p_k(\cdot;r)$ it follows that the set $\psi_k^{-1}(q;r,t)$ is a subset of this 2-dimensional space. Let 
\begin{equation}\label{xIsPreimage}
x\in \psi_k^{-1}(q;r,t).
\end{equation}
 Consider a 2-dimensional Cartesian coordinate system in it such that the point $o$ is on the ordinate axis and the point $x$ and the point that is symmetric to the point $x_k^*$ with respect to $o$ are on the abscissa axis (to determine the coordinate system uniquely we can additionally require the point $x_k^*$ to be in the upper half plane with respect to the abscissa axis). Let $x_k^*=(x_*, 2y_*)$, $q = (x_q,y_q)$, and $x = (x_1,0)$ in this coordinate system. Then $o=(0, y_*)$. By~\eqref{p_kDefinition} and the definition of the function $\psi_k$ we have that for arbitrary point $z = (x_z,y_z)$
 \begin{multline*}
 \psi_k(z;r,t) = (1-t)\frac {r(x_z+x_*, y_z+2y_*) + (0, y_* |x_k^*-z|_2)}{2r + |x_k^*-z|_2} + t(x_z,y_z) 
\\=
  \left((1-t)\frac{r(x_z+x_*)}{2r + |x_k^*-z|_2} + tx_z, (1-t) y_* + (1-t)\frac{ry_z}{2r + |x_k^*-z|_2}+ ty_z\right).
  \end{multline*}
  
From~\eqref{xIsPreimage} we obtain that $y_q = (1-t)y_*$; hence $\psi_k^{-1}(q;r,t)$ is a subset of the abscissa axis, and for all $z=(x_z,0)\in \psi_k^{-1}(q;r,t)$ we have $$(1-t)\frac{r(x_z+x_*)}{2r + \sqrt{(x_*-x_z)^2 + 4y^2_*}} + tx_z = x_q.$$

All the solutions of the latter equation are also solutions of the equation $$\left[(1-t)r(x_z+x_*) - 2r(x_q-tx_z)\right]^2 = (x_q-tx_z)^2  ((x_*-x_z)^2 + 4y^2_*),$$
which is an equation of not more than degree $4$ with respect to the unknown $x_z$; hence the set $\psi_k^{-1}(q;r,t)$ contains at most $4$ points. 

Finally, we construct the required partition of the set $U_k$. Consider the function $\psi_k(\cdot;r,t)$ as a function defined on the closure $\overline{U_k}$ of $U_k$. It is easy to see that $\psi_k$ is continuous, and hence by~\cite[Theorem~6.9.7]{Bogachev2} there exists a Borel set $B\subset \overline{U_k}$ such that $\psi_k(B;r,t) = \psi_k(\overline{U_k};r,t)$ and $\psi_k(\cdot ;r,t)$ is injective on $B$. 
Set $U_k^1 = B\cap U_k$. The set $\psi_k^{-1}(q;r,t)\setminus U_k^1$ consists of at most $3$ points for any $q\in \psi_k(U_k;r,t)$.
 Repeating the same  arguments we obtain measurable sets $U_k^2\subset U_k\setminus U_k^1$ and $U_k^3\subset U_k\setminus (U_k^1\cup U_k^2)$ such that $\psi_k$ is injective on each of them, and is injective on the measaruble set $U_k^4:= U_k\setminus (U_k^1\cup U_k^2\cup U_k^3)$.
The lemma is proved.
 \end{proof}

Below for a function $\phi\colon \Rd\to\Rd$ by $\frac{D\phi}{Dx}$ we denote its Jacobian matrix, $J_\phi(x):=\det\frac{D\phi}{Dx}(x)$ and $I$ is identity matrix. For a point or a vector $y\in\Rd$ and $s\in\{1,\dots, d\}$ by $y^s$ we denote its $s$-th coordinate.

\begin{lemma}\label{pr::pJacobian}
Assume a point $x\in U_k$ is given. Let $m$ be the middle of the segment $xx_k^*$, $\overline{\Delta x}=\frac {\overline{x_k^*x}}{|x-x_k^*|_2}$. Then
$$\frac {Dp_k(\cdot;r)}{Dx}(x) = \frac {r}{2r+|x-x_k^*|_2}I + \frac{2r}{(2r + |x-x^*_k|_2)^2}\cdot \left\|\overline{mo}^{\,i}\cdot\overline{\triangle x}^j\right\|_{i,j=1}^d.$$
\end{lemma}
\begin{proof} 
Let $x=(x^1,\dots, x^d)$, $x_k^* = (x_k^1,\dots,x_k^d)$, $o = (o^1,\dots, o^d)$. Let also $\delta_{ij} = 1$ if $i=j$, and  $\delta_{ij} = 0$ if $i\neq j$. Equality~\eqref{p_kDefinition} can be rewritten in the coordinates
$$p_k(x;r):=\frac{r\cdot(x^1+x_k^1,\dots, x^d+x_k^d)+ \sqrt{\sum_{s=1}^d(x^s-x_k^s)^2}\cdot (o^1,\dots, o^d)}
{\sqrt{\sum_{s=1}^d(x^s-x_k^s)^2} + 2r}.$$
\begin{gather*}
\frac{\partial p_k(x;r)^i}{\partial x^j} 
= 
\frac{(r\delta_{ij} + \frac{x^j-x^j_k}{|x-x^*_k|_2}o^i)\cdot (2r + |x-x^*_k|_2) - (r(x^i+x^i_k) + |x-x^*_k|_2o^i)\frac{x^j-x^j_k}{|x-x^*_k|_2} }
{(2r + |x-x^*_k|_2)^2}
\\=
\frac{r(2r+|x-x^*_k|_2)\delta_{ij}+2r\frac{x^j-x^j_k}{|x-x^*_k|_2}o^i-r\frac{x^j-x^j_k}{|x-x^*_k|_2}(x^i+x^i_k)}
{(2r + |x-x^*_k|_2)^2}
\\=
\frac {r}{2r + |x-x^*_k|_2}\delta_{ij} 
+\frac{2r}{(2r + |x-x^*_k|_2)^2}\cdot\frac{x^j-x^j_k}{|x-x^*_k|_2}\cdot\left(o^i - \frac{x^i+x^i_k}{2}\right).
\end{gather*}
The lemma is proved.  
\end{proof}

From Lemmas~\ref{pr::pJacobian} and~\ref{l::det} we immediately get the following lemma.
\begin{lemma}\label{pr::phiJacobian}
Let the function $\varphi_k(x;r,t)$ be defined in~\eqref{psiDefinition}. Then for all $t\in[0,1]$
$$J_{\varphi_k(\cdot;r,t)}(x) = t^d\frac{r^d}{(2r + |x-x^*_k|_2)^{d+1}}\left(2r+|x-x_k^*|_2 + 2 \left(\overline{mo}, \overline{\triangle x}\right)\right).$$
\end{lemma}

\begin{lemma}\label{pr::psiJacobian}
Let the function $\psi_k(x;r,t)$ be defined in~\eqref{psiDefinition}. Then for all $t\in[0,1]$
\begin{multline*}
J_{\psi_k(\cdot;r,t)}(x)
\\ = 
\left(\frac{r(1-t)}{2r + |x-x^*_k|_2} + t\right)^{d-1}
 \cdot
\left(\frac {2r^2 + 3r|x-x_k^*|_2 + |x-x_k^*|_2^2 - 2r\left(\overline{mo}, \overline{\triangle x}\right)}
{(2r+|x-x_k^*|_2)^2}(t-1) + 1\right).
\end{multline*}
\end{lemma}
\begin{proof}
To prove this lemma it is sufficient to apply Lemmas~\ref{pr::pJacobian} and~\ref{l::det} and notice that 
\begin{multline*}
\frac{r(1-t)}{2r + |x-x^*_k|_2} + t + \frac{2r(1-t)}{(2r + |x-x^*_k|_2)^2}\left(\overline{mo}, \overline{\triangle x}\right)
\\=
\frac {2r^2 + 3r|x-x_k^*|_2 + |x-x_k^*|_2^2 - 2r\left(\overline{mo}, \overline{\triangle x}\right)}
{(2r+|x-x_k^*|_2)^2}(t-1) + 1.
\end{multline*}
 \end{proof}
The next lemma justifies possibility to make a substitution  $y = \psi_k(x;r,t)$ in the integrals considered below.
\begin{lemma}
    For each fixed $0<r<R$ the function $\psi_k(x;r,t)$, $x\in U_k$, $t\in [0,1]$ is continuous on $U_k\times [0,1]$. The set $\Theta:=\{(x,t)\in U_k\times [0,1]\colon J_{\psi_k(\cdot;r,t)}(x) = 0\}$ has measure zero in $\RR^{d+1}$.
\end{lemma}
\begin{proof}
    Continuity of $\psi_k$ follows from the definition. The set $\Theta$ is a piece of the plot of the function
$$
t(x) = 1-\frac{(2r+|x-x_k^*|_2)^2} {2r^2 + 3r|x-x_k^*|_2 + |x-x_k^*|_2^2 - 2r\left(\overline{mo}, \overline{\triangle x}\right)}, x\in U_k
$$
on which $t(x)\in [0,1]$. There exists $\varepsilon>0$ such that 
$$
    \Theta\cap \left\{(x,t)\in U_k\times [0,1]\colon \left|2r^2 + 3r|x-x_k^*|_2 + |x-x_k^*|_2^2 - 2r\left(\overline{mo}, \overline{\triangle x}\right)\right| < \varepsilon\right\} = \emptyset.
$$
Thus on the set $\{x\in U_k\colon t(x)\in [0,1]\}$ the function $t$ is uniformly continuous, and hence its plot has zero measure.
\end{proof}

In the next two paragraphs we use the following notation.
The symbol $C$ stays for a positive number that does not depend on $n$. This number may be different in left and right parts of equality or inequality.
\subsection{Estimate for $\sum_{k=1}^{|S(n)|}\int_{T_k}|f(x) - f(p_k(x;r))|dx$}
\begin{lemma}\label{l::xpIntegralEstimate}
Assume $T_k\subset U_k$, $k=1,\dots, |S(n)|$, are measurable sets and 
\begin{equation}\label{diam<R}
\max_{k=1,\dots, |S(n)|}\diam U_k\leq r\leq R.
\end{equation}
 Then there exists a number $C$ that does not depend on $n$ and such that for all $f\in W^{1,p}(Q)$
$$\sum_{k=1}^{|S(n)|}\int_{T_k}|f(x) - f(p_k(x;r))|dx \leq Ch_n\left(\mes \bigcup_{k=1}^{|S(n)|} T_k\right)^{\frac 1 {p'}}  \|\,|\nabla f|_1\,\|_{L_p(Q)}.$$
\end{lemma}
\begin{proof}
Using Lemma~\ref{l::LagrangeTh} and Lemma~\ref{pr::smallDistance} we obtain
\begin{multline}\label{lxpk.1}
\sum_{k=1}^{|S(n)|}\int_{T_k}|f(x) - f(p_k(x;r))|dx 
\\ \leq 
\sum_{k=1}^{|S(n)|}\int_0^1\int_{T_k}|p_k(x;r) - x|_\infty |\nabla f((1-t)p_k(x;r) + tx)|_1dxdt
\\\leq 
Ch_n\sum_{k=1}^{|S(n)|}\int_0^1\int_{T_k} |\nabla f(\psi_k(x;r,t))|_1dxdt,
\end{multline}
where the functions $\psi_k(\cdot;r,t)$ are defined in~\eqref{psiDefinition}. Due to Lemma~\ref{pr::phiPsiInjectiveness} for each $t\in[0,1]$ and $k=1,\dots, |S(n)|$ there exists a partition $T_k = \bigcup_{i=1}^4 T_k^i(t)$ such that the function $\psi_k(\cdot;r,t)$ is injective on each of $T_k^i(t)$, $i=1,\dots, 4$. Hence for all $k=1,\dots, |S(n)|$
\begin{equation}\label{lxpk.2}
\int_0^1\int_{T_k} |\nabla f(\psi_k(x;r,t))|_1dxdt
 =
\sum_{i=1}^4
\int_0^1\int_{T_k^i(t)} |\nabla f(\psi_k(x;r,t))|_1dxdt,
\end{equation}
and for each $i=1,\dots, 4$, making a substitution $y=\psi_k(x;r,t)$ in the internal integral, and applying Holder's inequality, we can write
\begin{multline}\label{xp_kIntegralEstimate}
\int_0^1\int_{T_k^i(t)} |\nabla f(\psi_k(x;r,t))|_1dxdt
= 
\int_0^1\int_{\psi_k(T_k^i(t);r,t)} |\nabla f(y)|_1|J_{\psi_k^{-1} (\cdot;r, t)}(y)|dydt 
\\ \leq
\left(\int_0^1\int_{\psi_k(T_k^i(t);r,t)} |\nabla f(y)|_1^pdydt\right)^{\frac 1 p}
\left(\int_0^1\int_{\psi_k(T_k^i(t);r,t)}|J_{\psi_k^{-1}(\cdot;r, t)}(y)|^{p'}dydt\right)^{\frac {1}{p'}}.
\end{multline}

Applying the inverse function theorem and making a substitution $x=\psi_k^{-1}(y;r, t)$, we obtain
\begin{multline}\label{psiJacobianDoubleIntegral}
\int_0^1\int_{\psi_k(T_k^i(t);r,t)}|J_{\psi_k^{-1}(\cdot;r, t)}(y)|^{p'}dydt 
= \int_0^1\int_{\psi_k(T_k^i(t);r,t)}|J_{\psi_k(\cdot;r, t)}(\psi_k^{-1}(y;r, t))|^{-p'}dydt
\\=
\int_0^1\int_{T_k^i(t)}|J_{\psi_k(\cdot;r, t)}(x)|^{1-p'}dxdt
\leq
\int_{T_k}\int_0^1|J_{\psi_k(\cdot;r, t)}(x)|^{1-p'}dtdx.
\end{multline}
Since $p>d\geq 2$ we have $\frac {1}{p'}=1-\frac 1 p > \frac 1 2$ and hence 
\begin{equation}\label{p'Estimate}
0\geq 1-p'>-1. 
\end{equation}
Note that for all $t\in[0,1]$ 
$$
\frac{r(1-t)}{2r + |x-x^*_k|_2} + t= \frac{r(1+t) + t |x-x^*_k|_2}{2r + |x-x^*_k|_2}> \frac{r}{2r + r}=\frac 1 3, 
$$
\begin{multline*}
-\frac{\diam  Q}{2r}
\leq 
\frac {2r^2 + 3r|x-x_k^*|_2 + |x-x_k^*|_2^2 - 2r\left(\overline{mo}, \overline{\triangle x}\right)}
{(2r+|x-x_k^*|_2)^2}
\\\leq 
\frac {2r^2 + 3r^2 + r^2 + 2r \cdot \diam  Q}{4r^2} = \frac 3 2 + \frac {\diam  Q}{2r}.    
\end{multline*}

Using Lemma~\ref{pr::psiJacobian}, the latter inequalities and~\eqref{p'Estimate} we obtain 
\begin{equation}\label{psiJacobianIntegral}
\int_0^1|J_{\psi_k(\cdot;r, t)}(x)|^{1-p'}dt
\leq 
3^{(d-1)(p'-1)}\sup_{A}\int_0^1(A(t-1)+1)^{1-p'}dt\leq C,
\end{equation}
where the supremum is taken over the segment $A\in \left[-\frac{\diam  Q}{2r}, \frac 3 2 + \frac {\diam  Q}{2r}\right]$ and the number $C$ is finite and depends on $\diam  Q$ and $r$ and does not depend on $h_n$. Using~\eqref{psiJacobianIntegral}, \eqref{psiJacobianDoubleIntegral},  \eqref{xp_kIntegralEstimate}, \eqref{lxpk.2} and \eqref{lxpk.1}, we 	receive 
\begin{multline}\label{lxpk.3}
\sum_{k=1}^{|S(n)|}\int_{T_k}|f(x) - f(p_k(x;r))|dx
\\ \leq 
Ch_n\sum_{i=1}^{4}\sum_{k=1}^{|S(n)|}
\left(\int_0^1\int_{\psi_k(T_k^i(t);r,t)} |\nabla f(y)|_1^pdydt\right)^{\frac 1 p}(\mes T_k)^\frac{1}{p'}
\\ \leq 
Ch_n \sum_{i=1}^{4}
\left(\int_0^1\sum_{k=1}^{|S(n)|}\int_{\psi_k(T_k^i(t);r,t)} |\nabla f(y)|_1^pdydt\right)^{\frac 1 p}\left(\sum_{k=1}^{|S(n)|}\mes T_k\right)^\frac{1}{p'}.
\end{multline}
Note that for all $t\in [0,1]$ and each $i = 1,\dots, 4$ the sets $T_k^i(t)$ and the functions $\psi_k(\cdot; r,t)$ satisfy the conditions of Lemma~\ref{l::integralsSum}, because for all $x\in T_k^i(t)$ due to Lemma~\ref{pr::smallDistance}
\begin{multline*}
|x_k^*-\psi_k(x; r,t)|_\infty = |x_k^* - (1-t)p_k(x;r) - tx|_\infty 
\leq 
t|x_k^*-x|_\infty 
+
 (1-t)|x_k^*-p_k(x;r)|_\infty
\\\leq
|x_k^*-x|_\infty + |x_k^*-p_k(x;r)|_\infty
\leq
C|x_k^*-x|_\infty
\end{multline*}
with some number $C$ independent of $n$ (and $t$). 
To finish the proof of the lemma it is sufficient to apply Lemma~\ref{l::integralsSum} to~\eqref{lxpk.3}. The lemma is proved.  \end{proof}

\subsection{Estimate for $\sum_{k=1}^{|S(n)|}\int_{T_k}|f(x_k^*) - f(p_k(x;r))|dx$}
\begin{lemma}\label{l::x^*pIntegralEstimate}
Assume~\eqref{diam<R} holds, and measurable sets $T_k\subset U_k$, $k=1,\dots, |S(n)|$, are such that there exists a number $c>0$ that does not depend on $n$ for which $|J_{\varphi_k(\cdot;r,1)}(x)|>c$ for all $x\in T_k$. Then there exists a number $C$ that 
does not depend on $n$, such that for all $f\in W^{1,p}(Q)$
$$\sum_{k=1}^{|S(n)|}\int_{T_k}|f(x^*_k) - f(p_k(x;r))|dx \leq Ch_n\left(\mes \bigcup_{k=1}^{|S(n)|} T_k\right)^{\frac 1 {p'}}  \|\,|\nabla f|_1\,\|_{L_p(Q)}.$$
\end{lemma}
\begin{proof}
We may assume that each function $\varphi_k(\cdot;r,t)$  is injective on $T_k$, $k=1,\dots, |S(n)|$. Otherwise, due to Lemma~\ref{pr::phiPsiInjectiveness}, we can divide each of the sets $T_k$ into four subsets, so that the functions $\varphi_k(\cdot;r,t)$  are injective on each of the subsets, and apply arguments below to each of the subsets.
Using Lemma~\ref{l::LagrangeTh} we obtain
\begin{gather*}
\sum_{k=1}^{|S(n)|}\int_{T_k}|f(p_k(x;r))-f(x_k^*) |dx 
 \leq 
\sum_{k=1}^{|S(n)|}\int_0^1\int_{T_k}|p_k(x;r) - x_k^*|_\infty |\nabla f(\varphi_k(x;r,t))|_1dxdt
\\ \leq 
\sum_{k=1}^{|S(n)|}\int_0^1
\left(\int_{T_k}|p_k(x;r) - x_k^*|^{p'}_\infty dx\right)^{\frac 1 {p'}}
\left(\int_{T_k} |\nabla f(\varphi_k(x;r,t))|^p_1dx\right)^{\frac 1 {p}}dt
\\ \text{(by Lemma~\ref{pr::smallDistance})}
 \leq
C h_n\int_0^1\sum_{k=1}^{|S(n)|}(\mes T_k)^{\frac 1 {p'}}\left(\int_{T_k} |\nabla f(\varphi_k(x;r,t))|^p_1dx\right)^{\frac 1 {p}}dt
\\ \leq
C h_n\left(\sum_{k=1}^{|S(n)|}\mes T_k\right)^{\frac 1 {p'}}
\int_0^1\left(\sum_{k=1}^{|S(n)|}\int_{T_k} |\nabla f(\varphi_k(x;r,t))|^p_1dx\right)^{\frac 1 {p}}dt
\\ =
C h_n\left(\sum_{k=1}^{|S(n)|}\mes T_k\right)^{\frac 1 {p'}}\int_0^1\left(\sum_{k=1}^{|S(n)|}\int_{\varphi_k(T_k;r,t)} |\nabla f(y)|^p_1|J_{\varphi_k^{-1}(\cdot;r, t)}(y)|dy\right)^{\frac 1 {p}}dt
\\ =
C h_n\left(\sum_{k=1}^{|S(n)|}\mes T_k\right)^{\frac 1 {p'}}\int_0^1\left(
\sum_{k=1}^{|S(n)|}
\int_{\varphi_k(T_k;r,t)} \frac{|\nabla f(y)|^p_1dy}{|J_{\varphi_k(\cdot;r, t)}(\varphi_k^{-1}(y;r, t))|}\right)^{\frac 1 {p}}dt.
\end{gather*}
Due to Lemma~\ref{pr::phiJacobian} for all $k=1,\dots, |S(n)|$, $t\in[0,1]$ and $x\in T_k$, $J_{\varphi_k(\cdot;r,t)}(x)= t^d J_{\varphi_k(\cdot;r,1)}(x).$
Hence, using the conditions of the lemma we receive
\begin{gather*}
\sum_{k=1}^{|S(n)|}\int_{T_k}|f(p_k(x;r))-f(x_k^*) |dx 
\\ \leq 
C h_n\left(\sum_{k=1}^{|S(n)|}\mes T_k\right)^{\frac 1 {p'}}\int_0^1t^{-\frac d p}\left(\sum_{k=1}^{|S(n)|}
\int_{\varphi_k(T_k;r,t)} \frac{|\nabla f(y)|^p_1dy}{|J_{\varphi_k(\cdot;r, 1)}(\varphi_k^{-1}(y;r, t))|}\right)^{\frac 1 {p}}dt
\\ \leq 
C h_n\left(\sum_{k=1}^{|S(n)|}\mes T_k\right)^{\frac 1 {p'}}\int_0^1t^{-\frac d p}\left(\sum_{k=1}^{|S(n)|}
\int_{\varphi_k(T_k;r,t)} \frac{|\nabla f(y)|^p_1dy}{c}\right)^{\frac 1 {p}}dt.
\end{gather*}
Note that the sets $T_k$ and the functions $\varphi_k(\cdot;r,t)$ satisfy the conditions of Lemma~\ref{l::integralsSum}, since by Lemma~\ref{pr::smallDistance} for all $k=1,\dots, |S(n)|$, $t\in[0,1]$, and $x\in T_k$ 
$$|x_k^* - \varphi_k(x;r,t)|_\infty = t|p_k(x;r) - x_k^*|_\infty \leq |p_k(x;r) - x_k^*|_\infty \leq C|x-x_k^*|_\infty.$$
To finish the proof of the lemma it is sufficient to apply Lemma~\ref{l::integralsSum} and recall that $p>d$, and hence the integral $\int_0^1t^{-\frac d p}dt$ converges. The lemma is proved.  

\end{proof}

\subsection{Proof of Lemma~\ref{l::th5RemainderEstimate}}
\begin{proof}
For each $k=1,\ldots, |S(n)|$, we divide the set $U_k$ into three subsets 
$W_k^1:=\{x\in U_k\colon \left(\overline{mo}, \overline{\triangle x}\right) < -2R\}$, 
$W_k^2:=\{x\in U_k\colon \left(\overline{mo}, \overline{\triangle x}\right)\in [-2R, -\frac R 2]\}$ and $W_k^3:=\{x\in U_k\colon \left(\overline{mo}, \overline{\triangle x}\right) > -\frac R 2\}$.

Let $n$ be so big that for all $k=1, \dots, |S(n)|$ and $x\in U_k$, $|x-x^*_k|_2 <\frac R 8$.
For all $x\in W_k^1$ we have
$2R+|x-x_k^*|_2+ 2 \left(\overline{mo}, \overline{\triangle x}\right)<2R+\frac R 8-4R<-R,$
for all $x\in W_k^3$ we have
$2R+|x-x_k^*|_2+ 2 \left(\overline{mo}, \overline{\triangle x}\right)>2R-R=R,$
and hence by Lemma~\ref{pr::phiJacobian} for all $x\in W_k^1\cup W_k^3$ we have 
$$
\left|J_{\varphi_k(\cdot;R,1)}(x)\right|= \frac{R^d}{(2R + |x-x^*_k|_2)^{d+1}}\left|2R+|x-x_k^*|_2 + 2 \left(\overline{mo}, \overline{\triangle x}\right)\right|
 >
\frac {R^d}{(3R)^{d+1}}R=\frac {1}{3^{d+1}}.
$$

For all $x\in W_k^2$ we have
$$2\cdot \frac R 8+|x-x_k^*|_2+ 2 \left(\overline{mo}, \overline{\triangle x}\right)<\frac {2R} 8 +\frac R 8-R<-\frac R 2,$$
and hence by Lemma~\ref{pr::phiJacobian} for all $x\in W_k^2$ we have 
\begin{multline*}
\left|J_{\varphi_k\left(\cdot;\frac R 8,1\right)}(x)\right|= \frac{\left(\frac R 8\right)^d}{\left(\frac {2R} 8 + |x-x^*_k|_2\right)^{d+1}}\left|\frac {2R} 8+|x-x_k^*|_2 + 2 \left(\overline{mo}, \overline{\triangle x}\right)\right|
\\ >
\frac {8R^d}{(2R + R)^{d+1}}\cdot\frac R 2=\frac {4}{3^{d+1}}.
\end{multline*}

Finally let's return to the proof of~\eqref{th6.remainderEstimate}.
\begin{gather*}  
\sum_{k=1}^{|S(n)|}\int_{U_k} \left|f(x)-f(x_k^*)\right|dx
 =
\sum_{k=1}^{|S(n)|}\int_{W_k^1\cup W_k^3} \left|f(x)-f(x_k^*)\right|dx + \sum_{k=1}^{|S(n)|}\int_{W_k^2} \left|f(x)-f(x_k^*)\right|dx
\\ \leq 
\sum_{k=1}^{|S(n)|}\int_{W_k^1\cup W_k^3} \left|f(x)-f(p_k(x; R))\right|dx + 
\sum_{k=1}^{|S(n)|}\int_{W_k^1\cup W_k^3} \left|f(p_k(x; R))-f(x_k^*)\right|dx 
\\ + 
\sum_{k=1}^{|S(n)|}\int_{W_k^2} \left|f(x)-f\left(p_k\left(x;\frac R 8\right)\right)\right|dx + 
\sum_{k=1}^{|S(n)|}\int_{W_k^2} \left|f\left(p_k\left(x;\frac R 8\right)\right)-f(x_k^*)\right|dx.
\end{gather*}

To prove~\eqref{th6.remainderEstimate} we need to apply Lemmas~\ref{l::xpIntegralEstimate}, \ref{l::x^*pIntegralEstimate} and note that $\mes \bigcup_{k=1}^{|S(n)|} U_k = o(1)$ as $n\to\infty$ due to Property~\ref{U_k} of Lemma~\ref{l::InformationalSetProps}. The lemma is proved.  \end{proof}

\bibliographystyle{unsrt}
\bibliography{bibliography}
\end{document}